\documentclass{elsart}
\input{xypic}
\usepackage{amssymb}
\usepackage{amsmath}
\usepackage{latexsym}

\newtheorem{theorem}{Theorem}[section]
\newtheorem{proposition}[theorem]{Proposition}
\newtheorem{corrolary}[theorem]{Corollary}
\newtheorem{lemma}[theorem]{Lemma}
\newtheorem{example}[theorem]{Example}
\newtheorem{remark}[theorem]{Remark}
\newtheorem{question}[theorem]{Question}
\newtheorem{definition}[theorem]{Definition}
\newtheorem{notation}[theorem]{Notation}

\newcommand{\Z}{{\bf Z}}
\newcommand{\Q}{{\bf Q}}

\newcommand{\Hom}{{\rm Hom}}
\newcommand{\Ext}{{\rm Ext}}
\newcommand{\coker}{{\rm Coker}\,}
\newcommand{\im}{{\rm Im}\,}
\newcommand{\Ker}{{\rm Ker}\,}
\newcommand{\colim}{{\rm colim}\,}

\newcommand{\da}{\downarrow}
\newcommand{\ra}{\rightarrow}
\newcommand{\too}{\longrightarrow}
\newcommand{\lra}{\longrightarrow}
\newcommand{\epi}{\mbox{$\to$\hspace{-0.35cm}$\to$}}
\newcommand{\mono}{\hookrightarrow}

\newcommand{\map}{{\rm map}}
\newcommand{\Zup}{\Z[1/p]}
\newcommand{\Zpinf}{\Z(p^{\infty})}
\newcommand{\ov}{\overline}
\newcommand{\Gab}{G_{\rm ab}}

\begin{document}
\begin{frontmatter}

\title
{A connection between cellularization for groups and spaces
via two-complexes}

\author[Almeria]{Jos\'e L. Rodr\'{\i}guez\thanksref{rodriguez}},
\ead{jlrodri@ual.es}
\author[Barcelona]{J\'er\^{o}me Scherer\corauthref{cor}\thanksref{scherer}}
\ead{jscherer@mat.uab.es} \corauth[cor]{Corresponding author}

\address[Almeria]{\'{A}rea de Geometr\'{\i}a y Topolog\'{\i}a, Universidad
de Almer\'\i a, E--04120 Almer\'\i a}

\address[Barcelona]{Departament de Matem\`atiques,
Universitat Aut\`onoma de Barcelona, E--08193 Bellaterra}

\thanks[rodriguez]{Partially supported by the
Spanish Ministry of Education and Science MEC-FEDER grant
MTM2004-03629.}
\thanks[scherer]{Partially supported by the
program Ram\'on y Cajal, MEC, Spain, and the MEC-FEDER grant
MTM2004-06686.}


\begin{abstract}
Let $M$ denote a two-dimensional Moore space (so $H_2(M; \Z) =
0$), with fundamental group $G$. The $M$-cellular spaces are those
one can build from $M$ by using wedges, push-outs, and telescopes
(and hence all pointed homotopy colimits). The question we address
here is to characterize the class of $M$-cellular spaces by means
of algebraic properties derived from the group~$G$. We show that
the cellular type of the fundamental group and homological
information does not suffice, and one is forced to study a certain
universal extension.
\end{abstract}


\end{frontmatter}



\section*{Introduction}
Every pointed space $X$ can be approximated by an $M$-cellular
complex $CW_M X$ by means of a map $CW_MX\to X$ which induces a
weak homotopy equivalence on pointed mapping spaces
$$
\map_*(M, CW_M X) \to \map_*(M, X).
$$
This result generalizes the well known $CW$-approximation theorem
of J.H.C. Whitehead (in such case $M=S^1$). Bousfield was the
first to construct such a functor in the homotopy category,
\cite[Corollary~7.5]{Bou77}. Dror Farjoun's \cite{Far96} treats
extensively $CW_M$ as a functor in the category of pointed spaces
and studies its general properties in the context of homotopical
localization, see also \cite{Bla97} and \cite{Cha96}. In this
paper we focus on the case when $M$ is a two-dimensional Moore
space, i.e. a $2$-complex with $H_2(M; \Z) = 0$. We study the
relationship between the class of $M$-cellular spaces (those for
which $X \simeq CW_M X$) and the group theoretical properties of
$G = \pi_1 M$.

Let us denote by $J$ the set of primes $p$ for which $G_{\rm ab}$
is uniquely $p$-divisible. Define $R=\Z_{(J)}$, the integers
localized at $J$, if $G_{\rm ab}$ is torsion, and $R=\oplus_{p \in
J} \Z/p$ otherwise. Since $M$ itself is $HR$-acyclic, so is any
$M$-cellular space. We also noticed in \cite{RS99} that the
fundamental group of an $M$-cellular space is always a
$G$-cellular group. For a group $G$, the class of $G$-cellular
groups is, similarly as for spaces, the smallest class of groups
containing $G$ and closed under colimits. The group theoretical
cellularization has been first studied in \cite{RS99} and appears
in \cite{Ramon} in relation with $B\Z/p$-cellularization.
Recently, a purely algebraical approach by Farjoun, Dwyer,
G{\"o}bel, Shelah and Segev yielded important classification
results (\cite{FGS06-1}, \cite{FGS06-2}, \cite{DF05}).

The two observations about cellular spaces leads to a first and
naive guess for the characterization of $M$-cellular spaces:

\noindent {\bf Question.} Is a space $X$ $M$-cellular if and only
if $\pi_1 X$ is $G$-cellular and $X$ is $HR$-acyclic?

\medskip

The answer is yes when $M$ is the classical Moore space $M(\Z/p^n,
1)$, cofiber of the degree $p^n$ map on $S^1$, \cite{RS99},
extending results previously obtained in \cite{Bla97} and
\cite{CHR96} for Moore spaces $M(\Z/p^n, m)$ with $m\geq 2$. The
answer is actually positive for a larger class of Moore spaces:

\medskip

\noindent
{\bf Theorems \ref{divisiblecellular} and
\ref{finitetorsioncellular}.}
{\it Let $M$ be a two-dimensional Moore space, whose fundamental
group $G$ is either finite abelian or a subring of $\Q$. Then a
space $X$ is $M$-cellular if and only if $\pi_1 X$ is $G$-cellular
and $X$ is $HR$-acyclic.}

\medskip

In general the class of cellular spaces has not such a limpid
characterization. There is an extra condition on the second
homology group, related to the property of being ``quasi
$G$-radical", see Theorem~\ref{general-description}. However, the
characterization often holds for simply connected spaces.

\medskip

\noindent
{\bf Theorems \ref{rankonecellular} and \ref{torsioncellular}.}
{\it Let $M$ be a two-dimensional Moore space whose fundamental
group is either torsion abelian or a subgroup of $\Q$. Then a
simply connected space $X$ is $M$-cellular if and only if $X$ is
$HR$-acyclic.}

\medskip

Let us finally explain why the answer to the question is no in
general. Let $G=\Zup * \Z/p$ and $M=M(\Zup,1)\vee M(\Z/p,1)$. Then
$\Gab = \Zup \oplus \Z/p$ and the set of primes $J$ by which
$\Gab$ is uniquely divisible is the empty set. Therefore the
associated ring $R = \oplus_{q \in J} \Z/q$ is $0$, so the
condition $\tilde H_*(X;R)=0$ is always satisfied. However if the
naive characterization given above was true, all simply connected
spaces would be $M$-cellular. But this is false:

\medskip

\noindent
{\bf Theorem \ref{theorem counterexample}.}
{\it Let $G = \Zup * \Z/p$, and $M$ be the Moore space $M(\Z[1/p],
1) \vee M(\Z/p, 1)$, so the associated ring $R$ is $0$. The space
$K(\Z, 2)$ is then $HR$-acyclic, its fundamental group is $\Z[1/p]
* \Z/p$-cellular, but $K(\Z, 2)$ is not $M$-cellular.}

\medskip

Even though this result provides a negative answer to the question
we asked above, quite a few problems remain open. We ask in
particular whether the naive characterization holds for Moore
spaces with abelian fundamental group. This brings us to a
question which can be considered as purely group theoretical, see
Question~\ref{listofMoorespaces}. If $M(G, 1)$ is a
two-dimensional Moore space with abelian fundamental group $G$, is
$G$ a quotient of a subgroup of $\Q$?


\section{A general study of $M(G,1)$-cellular spaces}

In this section $M$ is a two-dimensional Moore space $M(G,1)$,
where $G \cong \pi_1 M(G,1)$. The objective is to obtain a
characterization of $M(G, 1)$-cellular spaces.

\subsection{Chach\'olski's approach to cellularization}
The fundamental theorem on which our work on cellularization
relies is that of Chach\'olski stating that $CW_M X$ can be
obtained as a certain homotopy fibre. Recall that $P_{M}$ denotes
the $M$-nullification functor, see \cite{Far96}. In particular
connected spaces $X$ such that $P_M X \simeq X$ are called
$M$-\emph{null} and are characterized by the fact that $\map_*(M,
X) \simeq *$.

\begin{theorem}
\label{Chacholski}{\rm \cite[Theorem 20.3]{Cha96}}
Let $X\to X'$ be a map which induces the trivial map $[M,X]\to [M,
X']$ and assume that its homotopy fiber is $M$-cellular. Then
$CW_MX$ is the homotopy fibre of the composite map $X \to X' \to
P_{\Sigma M} X'$. \hfill\qed
\end{theorem}

As explained in \cite[Theorem 20.5]{Cha96}, one way to construct a
map $X \rightarrow X'$ having the desired properties is the
following: Take a wedge of as many copies of $M$ as there are
homotopy classes of maps $M \ra X$ and consider the cofibration
sequence
\begin{equation}
\label{xprime} \bigvee_{[M, X]} M \stackrel{\varepsilon}{\ra} X
\ra X'=P^1_M X,
\end{equation}
where the map $\varepsilon$ is the evaluation map. The homotopy
fiber of $X \to P^1_M X$ is then $M$-cellular and $CW_M X$ is
obtained as the homotopy fiber of a mixing process between the
$M$-nullification and the $\Sigma M$-nullification. Indeed $P^1_M
X$ can be seen as the first step of the construction of $P_M X$
and $CW_M X \simeq Fib(X \ra P^1_M X \ra P_{\Sigma M} P^1_M X)$.
In particular, a space $X$ is $M$-cellular if and only if the
space $P^1_M X$ is $\Sigma M$-acyclic, i.e. $P_{\Sigma M} P^1_M X
\simeq *$. We show next that one can sometimes make do with less.

\begin{corrolary}\label{wojtekbetter}
Let $M$ be any space and assume there exists a subset $I$ of $[M,
X]$ such that the homotopy cofiber $C_\epsilon$ of the evaluation
map $\epsilon: \bigvee_{I} M \ra X$ is already $\Sigma M$-acyclic.
Then $X$ is $M$-cellular.
\end{corrolary}

\begin{pf}
We have to show that $P^1_M X$ is also $\Sigma M$-acyclic. For
this purpose we use the cofibration $\bigvee_{I} M \mono
\bigvee_{[M, X]} M \ra \bigvee_{I'} M$, where $I'$ is the
complement of $I$ in $[M, X]$. This induces a cofibration sequence
$C_\epsilon \ra P^1_M X \ra \bigvee_{I'} \Sigma M$. But since
$C_\epsilon$ is killed by $\Sigma M$, it is easy to check that so
is $P^1_M X$. \hfill{\qed}
\end{pf}

\subsection{The characterization of $M(G,1)$-cellular spaces}
{}From now on $M$ is a two-dimensional Moore space $M(G,1)$. The
homology exact sequence associated to the cofibration
(\ref{xprime}) has then the following simple form:
\begin{equation}
\label{exact-sequence}
 0\to H_2 X \to E \to \bigoplus_{[M,X]}  \Gab \to H_1X \to 0
\end{equation}

\begin{notation} \label{notation}
{\rm Let $G$ be a group. As for example in \cite[Section
7]{Bou97}, $J$ denotes the set of primes $p$ for which $G_{\rm
ab}$ is uniquely $p$-divisible, and $J'$ its complementary set of
primes. If $G_{\rm ab}$ is torsion, then we set $H=\oplus_{p \in
J'} \Z/p$ and $R=\Z_{(J)}$. In the other case, i.e., if $\Gab$
contains some torsion free element, $H=\Z[J^{-1}]$ and
$R=\oplus_{p \in J} \Z/p$.}
\end{notation}

\begin{lemma} \label{casacuberta}
The map $[M, X] \ra \Hom(G, \pi_1 X)$ is an epimorphism for
any~$X$. \hfill\qed
\end{lemma}

Recall that the $G$-{\em radical} $T_G N$ of a group $N$ is the
smallest subgroup of $N$ such that $\Hom(G, N/T_G N) = 0$,
\cite{Bou77}. When $N = T_G N$, one says that $N$ is $G$-radical.

\begin{theorem} \label{general-description}
A space $X$ is $M$-cellular if and only if the following three
conditions are satisfied: \begin{itemize}
\item[(1)] the fundamental group $\pi_1X$ is
$G$-cellular,

\item[(2)] the space $X$ is $HR$-acyclic,

\item[(3)] the abelian group $E$ in
(\ref{exact-sequence}) is $\Gab$-radical.

\end{itemize}
\end{theorem}

\begin{pf}
According to Proposition 5.3 in \cite{RS99} the space $P_{\Sigma
M} P^1_M X$ is contractible if and only if $P^1_M X$ is
1-connected, $\pi_2(P^1_M X)$ is $\Gab$-radical and $\tilde
H_*(P^1_M X;R)=0$. Since $M$ is $HR$-acyclic, the cofibration
(\ref{xprime}) tells us that $X$ is $HR$-acyclic if and only if
$P^1_M X$ is $HR$-acyclic. The associated homology exact sequence
in low dimensions is (\ref{exact-sequence}), where $E= H_2(P^1_M
X)\cong \pi_2(P^1_M X)$ and $P^1_M X$ is 1-connected by the above
lemma. The conclusion of the theorem is then clear. \hfill{\qed}
\end{pf}

\begin{corrolary}\label{1-connected-description}
A 1-connected space $X$ is $M$-cellular if and only if $X$ is
$HR$-acyclic and the abelian group $E$ in (\ref{exact-sequence})
is $\Gab$-radical. \hfill{\qed}
\end{corrolary}

\subsection{Quasi-radical groups}
Let us remark here that in the situation of
Corollary~\ref{1-connected-description}, $E\cong \colim_\alpha
(H_2X \mono E_\alpha)$ where $\alpha$ runs over all extensions of
$\Gab$ by $H_2X$. This is a ``universal" extension in the sense
that any extension $H_2X \mono E_\alpha \epi \Gab$ can be obtained
as the pull-back of it along the inclusion $i_\alpha: \Gab \mono
\oplus \Gab$.

\begin{definition} \label{quasiradical}
{\rm Let $G$ be any group, and let $H$ be the group associated to
$\Gab$ as in Notation~\ref{notation}. An $H$-radical abelian group
$A$ is called {\it quasi $G$-radical} if the universal extension
$$
0\to A \to E \too \bigoplus_{\Ext(\Gab,A) } \Gab \to 0
$$
is $\Gab$-radical.}
\end{definition}

\begin{remark} \label{strictinclusion}
{\rm There is a strict inclusion \{$\Gab$-radical groups\}
$\subseteq$ \{quasi $G$-radical groups\}. For $G= \Z_{p^\infty}$
for example, the group $H$ is $\Z/p$. Is it itself
quasi-$G$-radical, but not $G$-radical.}
\end{remark}

Our terminology is the analogue of Mislin and Peschke's notion of
$h2$-\emph{perfect} groups given in \cite[Section~2]{MP01}. In
fact the two settings intersect in the most simplest cases: When
$G = \Z/p$, an abelian group is quasi $\Z/p$-radical if and only
if it is $H\Z[1/p]\, 2$-perfect (compare \cite[Example~5.1]{MP01}
with Theorem~\ref{finitetorsioncellular}) and likewise, when $G =
\Z[1/p]$, choose $H\Z/p$, ordinary homology with coefficients in
$\Z/p$ (compare \cite[Example~5.2]{MP01} with
Proposition~\ref{divisiblecellularization}).

\begin{proposition}
\label{propK(A, 2)}
Let $A$ be an abelian group such that $K(A,2)$ is $M$-acyclic.
Then $K(A,2)$ is $M$-cellular if and only if $A$ is quasi
$G$-radical. \hfill{\qed}
\end{proposition}

The study of the cellularization of Eilenberg-Mac Lane spaces will
be refined in Section~\ref{sec counter}.
\medskip

Since the integral homology groups of an $HR$-acyclic space are
$H$-radical (by the universal coefficients theorem), we obtain the
following reformulation of
Corollary~\ref{1-connected-description}:

\begin{corrolary}\label{1-connected-description2}
A simply connected space $X$ is $M$-cellular if and only if
\begin{itemize}
\item[{\rm (1)}] the group $\pi_2 X$ is quasi $G$-radical, and

\item[{\rm (2)}] $\pi_k X$ is $J'$-torsion for all $k$, if $\Gab$ is torsion,
or uniquely $J$-divisible otherwise. \hfill\qed
\end{itemize}
\end{corrolary}

\subsection{Consequences for the naive characterization}
It is difficult to know whether a group is quasi-radical at a
first glance. Therefore it is useful to have partial answers to
our question where some extra assumption on $X$ forces $\pi_2 X$
to be quasi-radical. All the cases where the naive
characterization holds are obtained by comparing $H$-radical
groups with quasi $G$-radical ones.

\begin{theorem}\label{H2-G-radical}
Let $X$ be a space such that $H_2X$ is $\Gab$-radical. Then $X$ is
$M$-cellular if and only if $\pi_1X$ is $G$-cellular and $X$ is
$HR$-acyclic.
\end{theorem}

\begin{pf}
Suppose that $X$ is $HR$-acyclic, the fundamental group $N$ is
$G$-cellular, and $H_2X$ is $\Gab$-radical. Let us consider a
little variation of the construction of $CW_MX$ (which in general
does not produces $CW_MX$ in the sense of
Theorem~\ref{Chacholski}). Instead of $P^1_M X$, we take
$C_\epsilon$ the homotopy cofibre of the evaluation map
$\bigvee_IM\to X$ where $I$ is the union of the maps in $[M,X]$
represented by non-trivial elements in $\Ext(\Gab, \pi_2X)$ and
one preimage under $[M,X] \epi \Hom(G,N)$ for each morphism $G \ra
N$. By Theorem~\ref{wojtekbetter}, $X$ is $M$-cellular if we show
that $C_\epsilon$ is $\Sigma M$-acyclic. The map $X\to K(N,1)$
classifying the universal cover of $X$ yields a diagram whose rows
are cofibrations
$$
\begin{array}{ccccc}
\bigvee_{I} M & \to  & X  & \to  & C_\epsilon\\
\da & & \da & & \da \\
\bigvee_{\Hom(G,N)}  M & \to  & K(N,1)   & \to  & K(N,1)'.\\
\end{array}
$$
Thus we get a commutative diagram in homology where the rows are exact
$$
\begin{array}{ccccccccc}
0 \to & H_2X &\to & H_2(C_\epsilon)  & \to  & \bigoplus_I \Gab &
\stackrel{h''}{\longrightarrow}   & H_1N &\to 0\\
& \da &&\da & & \da & & \| & \\
0 \to &H_2 N & \to & H_2 (K(N,1)')& \stackrel{g}{\lra}
&\bigoplus_{\Hom(G,N)}  \Gab  & \stackrel{h'}{\longrightarrow} & H_1N &\to 0
...\\
\end{array}
$$
The kernel of $h''$ is the direct sum of $\bigoplus_{\Ext(\Gab,
\pi_2X)} \Gab$ and the kernel of $h'$. Now, since $N$ is
$G$-cellular by assumption, the group $H_2(K(N,1)')$ is
$\Gab$-radical by \cite[Corollary 3.8]{RS99}, and thus $\im g=\Ker
h'$ is also $\Gab$-radical. Hence, $H_2(C_\epsilon)$ is
$\Gab$-radical and this proves that $C_\epsilon$ is  $\Sigma
M$-acyclic. \hfill{\qed}
\end{pf}

The condition of $H_2X$ being $\Gab$-radical is not necessary. For
instance, we will see in Theorem~\ref{pinfcellular} that
$X=M(\Z/p,2)$ is $M(\Zpinf,1)$-cellular, while $\Z/p$ is not
$\Zpinf$-radical. It is of course quasi $\Zpinf$-radical, compare
with Remark~\ref{strictinclusion}.

\begin{corrolary}
\label{corollary_Gradical}
Let $X$ be a space such that $\pi_2X$ is $\Gab$-radical. If
$\pi_1X$ is $G$-cellular and $X$ is $HR$-acyclic, then $X$ is
$M$-cellular.
\end{corrolary}

\begin{pf}
In the Hopf exact sequence $\pi_2X\to H_2X\epi H_2N$, both
$\pi_2X$ and $H_2N$ (by \cite[Corollary 2.8]{RS99}) are
$\Gab$-radical. Thus so is $H_2X$ and we conclude by
Theorem~\ref{H2-G-radical}. \hfill{\qed}
\end{pf}

In fact a slight variation in the proof of
Theorem~\ref{H2-G-radical} yields the following generalization of
the preceding corollary.

\begin{proposition}\label{pi2-quasi-G-radical}
Let $X$ be a space such that $\pi_2X$ is quasi-$G$-radical. Then
$X$ is $M$-cellular if and only if $\pi_1X$ is $G$-cellular and
$X$ is $HR$-acyclic. \hfill{\qed}
\end{proposition}


\section{Moore spaces with abelian fundamental group}

In this section we focus on two-dimensional Moore spaces $M(G, 1)$
with abelian fundamental group~$G$. This study will be used in the
last section to give counter-examples to the naive
characterization of $M(G, 1)$-cellular spaces. Let us give
concrete models for all known two-dimensional Moore spaces with
abelian fundamental group.

\subsection{Varadarajan's list}
We start with the following list of all abelian groups $G$ with
$H_2(G; \Z)=0$. These are exactly the groups for which a Moore
space $M(G,1)$ exists. It might however be a 3-complex.

\begin{theorem}\label{varadarajan}\cite[Theorem 2.6]{Var66}
Let $G$ be an abelian group, with $T$ the torsion subgroup of~$G$.
For any prime $p$, let $T(p)$ denote the $p$-primary component
of~$T$. Then there exists a Moore space $M(G, 1)$ if and only if
\begin{itemize}
 \item[(1)] The quotient $G/T$ is of rank at most 1 over $\Q$,
 \item[(2)] for all $p$, the group $T(p)$ is either divisible or the direct sum
 of a divisible group and a cyclic group,
 \item[(3)] for all $p$, we have $T(p) \otimes G/T = 0$. \hfill\qed
\end{itemize}

\end{theorem}

\subsection{The torsion free case}
Let $S$ be any subgroup of $\Q$, the ring of rationals. We can
assume that $1$ is in~$S$. These groups have been classified by
Baer and are determined by their type (see for example
\cite[Theorem 85.1]{Fu73}). The {\it type} of $S$ is the sequence
$(k_2, k_3, k_5, k_7, k_{11}, \cdots, k_p, \cdots)$ where $p$ runs
over the set $\mathcal P$ of all primes and $k_p$ is either a
natural integer or infinite. The number $k_p$ indicates that $1$
is divisible in $S$ by $p^{k_p}$, but not by $p^{k_p+1}$, unless
$k_p$ is infinite. For example the group of type $(0, 0, \cdots,
0, \infty, 0, \cdots)$ is $\Z[1/p]$.

We now exhibit a particular construction of a Moore space $M(S,
1)$, where $S$ is a subgroup of $\Q$ of type $(k_p)_{p \in
\mathcal P}$. Let us first order all primes and their powers up to
$p^{k_p}$ by increasing order and we denote this sequence by
$(m_1, m_2, m_3, \cdots)$. We define then $\alpha_n$ to be the
unique prime $p$ dividing $m_n$. In this way, the sequence
$(\alpha_1, \alpha_2, \alpha_3, \cdots)$ contains exactly $k_p$
times the prime $p$. Obviously the colimit of
$$
\Z \stackrel{\alpha_1}{\lra} \Z \stackrel{\alpha_2}{\lra} \Z
\stackrel{\alpha_3}{\lra} \Z \stackrel{\alpha_4}{\lra} \cdots
$$
is $S$. We can also realize this at the level of spaces by a
telescope of circles
$$
S^1 \stackrel{\alpha_1}{\lra} S^1 \stackrel{\alpha_2}{\lra} S^1
\stackrel{\alpha_3}{\lra} S^1 \stackrel{\alpha_4}{\lra} \cdots
$$
and its homotopy colimit is $M(S, 1)$. It has dimension 2 since it
is a telescope of one dimensional spaces.

\subsection{The torsion case}
Theorem~\ref{varadarajan} tells us that there exist very few Moore
spaces with abelian torsion fundamental group. Such a group has to
be either divisible, or a direct sum with a cyclic group. The
classical Moore spaces, homotopy cofiber of the $n$-th power map
on the circle, are Moore spaces $M(\Z/n, 1)$ and these are
two-dimensional. There exists also a two-dimensional Moore space
for $\Zpinf$, the Pr\"ufer group, cokernel of the canonical
inclusion of $\Z$ in $\Z[1/p]$. This is a particular case of the
following construction, which uses the Moore spaces $M(S, 1)$ we
defined above, where $S$ is any subgroup of $\Q$.

Let $S$ be a subgroup of $\Q$ of type $(k_p)_{p \in \mathcal P}$.
The set $J$ for $S$ consists in the primes $p$ for which $k_p$ is
infinite, and the cokernel of the inclusion $\Z \mono S$ is
isomorphic to $\oplus_{p \in J} \Zpinf \oplus \oplus_{q \in J'}
\Z/q^{k_q}$. Taking the corresponding element $\alpha$ in $\pi_1
M(S, 1)$, we get a cofibration sequence
$$
S^1 \stackrel{\alpha}{\lra} M(S, 1) \lra M(\bigoplus_{p \in J}
\Zpinf \oplus \bigoplus_{q \in J'} \Z/q^{k_q}, 1).
$$
This last space is again a 2-complex.

\begin{question}\label{listofMoorespaces}
{\rm Are there any other two-dimensional Moore spaces with abelian
fundamental group? In particular we do not know if there exists
one whose fundamental group is isomorphic to $\Z(p^\infty) \oplus
\Z(p^\infty)$.}
\end{question}

This question could be considered as purely group theoretical. The
existence of a Moore space with fundamental group $G$ is indeed
equivalent to the existence of a presentation $* \Z
\xrightarrow{\alpha} * \Z \rightarrow G$ such that the
abelianization $\alpha_{\rm ab}: \oplus \Z \rightarrow \oplus \Z$
is injective.

\subsection{Cellularity in the non-torsion case}
The set of primes $J$ associated to a subgroup $S < \Q$ consists
in those $p \in \mathcal P$ such that $k_p = \infty$. The group
$H$ is then $\Z[J^{-1}]$ and the ring $R = \oplus_J \Z/p$. Thus a
space $X$ is $HR$-acyclic if and only if its integral homology is
uniquely $J$-divisible.

\begin{lemma}\label{modnrankone}
Let $M$ be a two-dimensional Moore space with fundamental group a
subgroup $S$ of $\Q$, of type $(k_p)_{p \in \mathcal P}$ and let
$n$ be an integer such that $p | n$ only if $p \not\in J$. Then
$M(\Z/n, 1)$ is $M$-cellular.
\end{lemma}

\begin{pf}
For such an integer $n$, the cokernel of the multiplication by $n$
on $S$ is $\Z/n$. \hfill{\qed}
\end{pf}

\begin{theorem}\label{rankonecellular}
Let $M$ be a two-dimensional Moore space with fundamental group a
subgroup $S$ of $\Q$. Then a 1-connected space $X$ is $M$-cellular
if and only if $X$ is $HR$-acyclic.
\end{theorem}

\begin{pf}
We have to prove that $S^2[J^{-1}]=M(\Z[J^{-1}],2)$ is $M$-cellular.
This is sufficient as the class
of $S^2[J^{-1}]$-cellular spaces coincides with the class of 1-connected
$HR$-acyclic spaces.

Let us look at the inclusion $\Z[J^{-1}] \mono S$. This is induced by a
morphism
of telescopes as follows. Define
$$
\beta_k = \left\{
            \begin{array}{lll}
              \alpha_k & {\hbox {\rm if}} & \alpha_k \in J \\
              1 & {\hbox {\rm if}} & \alpha_k \not\in J \\
            \end{array}
           \right.
$$
and define then by induction on $k$ a sequence $\gamma_k$, by $\gamma_0 =
1$ and
$\gamma_k = \gamma_{k-1} \cdot \beta_k$. We have then a commutative diagram
$$
\begin{array}{ccccccc}
\Z & \stackrel{\alpha_1}{\lra}  & \Z & \stackrel{\alpha_2}{\lra}  & \Z &
\stackrel{\alpha_3}{\lra}  & \cdots\\
\da^{\gamma_0} & & \da^{\gamma_1} & & \da^{\gamma_2} \\
\Z & \stackrel{\beta_1}{\lra}  & \Z & \stackrel{\beta_2}{\lra}  & \Z &
\stackrel{\beta_3}{\lra}  & \cdots\\
\end{array}
$$
Replacing every copy of the integers by a copy of a circle $S^1$,
we get the map $S^2[J^{-1}] \ra M$ as a map between telescopes.
Its homotopy cofiber is thus the homotopy colimit of the homotopy
cofibers of $\gamma_k: S^1 \ra S^1$. These are Moore spaces of
type $M(\Z/\gamma_k, 1)$. Notice that no prime $p \in J$ divides
any $\gamma_k$, i.e. $\gamma_k$ satisfies the conditions of
Lemma~\ref{modnrankone}. Hence $M(\Z/\gamma_k, 1)$ is $M$-cellular
and their telescope as well: $M(\oplus_{p \not\in J} \Z/p^{k_p},
1)$ is $M$-cellular. The Puppe sequence tells us finally that we
have a cofibration
$$
M \lra M(\oplus_{p \not\in J} \Z/p^{k_p}, 1) \lra S^2[J^{-1}]
$$
where the first two spaces are $M$-cellular. Thus so is the third
and we are done. \hfill{\qed}
\end{pf}

\begin{example}\label{type1}
{\rm Let $S$ be the subgroup of $\Q$ of type $(1,1,1,1, \cdots)$,
that is, $S$ is the additive subgroup of $\Q$ generated by $1/p$
for all prime numbers $p$. The above theorem shows in this case
that $S^2$ is $M(S, 1)$-cellular. There is an isomorphism $\Ext
(S, \Z) \cong (\prod_p \Z/p)/\Z$ and there are two possibilities
for an extension $\Z \ra E \ra S$. Either the $S$-reduction of~$E$
(the quotient by its $S$-radical) is $E$, or it is $\Z$, depending
on whether the element $\alpha \in (\prod_p \Z/p)/\Z$ representing
the extension is torsion or not. Thus none of these extensions is
$S$-radical. However, if we take $F =\colim_{E \in \Ext (S,
\Z)}(\Z \mono E)$, then $F$ is $S$-radical (i.e. $\Z$ is
$S$-quasi-radical).}
\end{example}

\subsection{The case of subrings of $\Q$}
We have here a stronger result than Theorem~\ref{rankonecellular}
and can in fact completely determine the class of $M$-cellular
spaces.

\begin{theorem}\label{divisiblecellular}
Let $M = M(\Z[J^{-1}], 1)$ be a two-dimensional Moore space. Then
a space $X$ is $M$-cellular if and only if $\pi_1 X$ is
$\Z[J^{-1}]$-cellular and $X$ is $HR$-acyclic.
\end{theorem}

\begin{pf}
This result is a reflection of the fact that the group $H$
associated to $\Z[J^{-1}]$ as introduced in
Notation~\ref{notation} is $\Z[J^{-1}]$ itself. This means we can
directly apply Theorem~\ref{H2-G-radical} because $H_2(X; \Z)$ is
$H$-radical. \hfill{\qed}
\end{pf}

\begin{lemma}\label{divisibleclasses}
Let $G = \Z[J^{-1}]$ and $M = M(G, 1)$ be a two-dimensional Moore
space. The class of nilpotent $M$-cellular spaces coincides then
with that of nilpotent $M$-acyclic spaces.
\end{lemma}

\begin{pf}
An $M$-cellular space is always $M$-acyclic,
Theorem~\ref{Chacholski}. Let thus $X$ be an $M$-acyclic space. We
apply Chach\'olski's theorem to prove it is also $M$-cellular. The
cofiber $P^1_M X$ of the evaluation $\vee M \ra X$ is a simply
connected space because there is no difference between $G$-socular
and $G$-radical nilpotent groups (there is none for abelian groups
and use \cite[Proposition 7.4]{Bou97}). It is also $M$-acyclic and
hence $H_n(P^1_M X;\Z)$ has to be uniquely $J$-divisible for all
$n \geq 2$, i.e. $G$-radical. Therefore $P_{\Sigma M} P^1_M X
\simeq *$ and we are done. \hfill{\qed}
\end{pf}

\begin{proposition}\label{divisiblecellularization}
Let $G = \Z[J^{-1}]$ and $M = M(G, 1)$ be a two-dimensional Moore
space. Then, for any nilpotent space $X$ with $G$-radical
fundamental group, we have $CW_M X \simeq Fib(X \lra \prod_J
X^{\;\widehat{}}_p \,)$.
\end{proposition}

\begin{pf}
This is a direct consequence of the preceding lemma and
\cite[Theorem 4.4]{CR97}, which identifies $P_M X$ with $\prod_J
X^{\;\widehat{}}_p$ when $X$ is a nilpotent space with $G$-radical
fundamental group. \hfill{\qed}
\end{pf}

\begin{example}\label{S^2}
{\rm Let $G=\Q$. Then $CW_M S^2$ is a two-stage Postnikov space, with
$$
\pi_1 CW_M S^2 \cong \pi_2 CW_M S^2 \cong (\prod_{\mathcal P}
\Z^{\;\widehat{}}_p)/\Z \cong \oplus \Q.
$$
This gives in particular an example where the fundamental group of
the $M$-cellulari\-za\-tion drastically differs from the
$G$-cellularization of the fundamental group.}
\end{example}

\begin{example}\label{Zpinf}
{\rm Let $G=\Z[1/p]$. Then $CW_M K(\Zpinf, 1) \simeq
K(\Q^{\;\widehat{}}_p, 1)$ and the associated universal central
extension (as in \cite[Theorem 2.7]{RS99}) is the extension
$\Z^{\;\widehat{}}_p \ra \Q^{\;\widehat{}}_p \ra \Zpinf$. }
\end{example}

\subsection{Cellularity in the reduced torsion case}
The naive description of the class of $M(G, 1)$-cellular spaces
given in the introduction is an actual characterization when $G$
is any finite cyclic group. More generally we have:

\begin{theorem}\label{finitetorsioncellular}
Let $G= \oplus \Z/p^{k_p}$ be an abelian torsion group with no
divisible summand and $M=M(G,1)$ be a two-dimensional Moore space.
Then a space $X$ is $M$-cellular if and only if $\pi_1 X$ is
generated by elements of order $p^{l}$ for $l \leq k_p$ and $X$ is
$H\Z_{(J)}$-acyclic.
\end{theorem}

\begin{pf}
The hypothesis on the fundamental group tells us that the cofiber
$X'$ in Theorem~\ref{Chacholski} is simply connected. Moreover,
since $X$ is $H\Z_{(J)}$-acyclic, so is~$X'$. We conclude by
Bousfield's explicit computations in \cite[Theorem~7.5]{Bou97}
that $P_{\Sigma M} X'$ is contractible. \hfill{\qed}
\end{pf}

For cyclic groups of prime power order, we recover \cite[Theorem
6.2]{RS99}.

\subsection{Cellularity in the unreduced torsion case}
When $G = \Zpinf$, the same characterization holds in the simply
connected case.

\begin{theorem}\label{pinfcellular}
Let $M=M(\Zpinf, 1)$. Then a 1-connected space $X$ is $M$-cellular if and
only if $X$ is $H\Z[1/p]$-acyclic.
\end{theorem}

\begin{pf}
The space $M(\Z/p^k, 2)$ is the homotopy cofiber of the map $M \ra
M$ induced by the multiplication by $p^k$ on $\Zpinf$, so it is
$M$-cellular. Hence any simply connected $p$-torsion space is so
(see \cite{CHR96}). \hfill{\qed}
\end{pf}

We will need the following lemma to analyze the class of $M(G,
1)$-cellular spaces for a general abelian torsion group~$G$.

\begin{lemma}\label{directsum}
Let $A$ and $B$ be two abelian groups such that there exist a
two-dimensional Moore space $M$ of type $M(A \oplus B, 1)$. Then
both $M(A, 1)$ and $M(B, 1)$ are $M$-cellular.
\end{lemma}

\begin{pf}
Let $M \ra M(A, 1)$ be any map which induces the projection $A
\oplus B \ra A$ on the fundamental group. Its homotopy cofiber is
$M(B, 2)$, which is $\Sigma M$-acyclic. Thus $M(A, 1)$ is
$M$-cellular by Chach\'olski's theorem, or rather by its
Corollary~\ref{wojtekbetter}. \hfill{\qed}
\end{pf}

\begin{remark}\label{retract}
{\rm We point out that it is not clear whether $M(A, 1)$ is a
retract of $M(A \oplus B, 1)$, which would provide a direct proof
of the above lemma. One can of course construct maps $M(A, 1)
\rightarrow M(A \oplus B, 1)$ and $M(A \oplus B, 1) \rightarrow
M(A, 1)$ which induce the canonical inclusion and projection, but
the composite might fail to be the identity, compare with
Lemma~\ref{casacuberta}.}
\end{remark}

We can now prove the same result as Theorem \ref{pinfcellular} for
any $M(G, 1)$ with $G$ torsion and abelian. Together with
Theorem~\ref{rankonecellular}, it shows that the naive
characterization given in the introduction holds for all Moore
spaces constructed at the beginning of this section, at least in
the 1-connected case.

\begin{theorem}\label{torsioncellular}
Let $M=M(G, 1)$ be a two-dimensional Moore space with torsion
abelian fundamental group $G$. Then a 1-connected space $X$ is
$M$-cellular if and only if $X$ is $H\Z_{(J)}$-acyclic.
\end{theorem}

\begin{pf}
By Theorem \ref{varadarajan} we infer that $T(p)$ is either
cyclic, or the direct sum of a divisible group and a cyclic one,
for all $p \in J$. Because $\Ext(T(p), T(q))$ is zero, our group
$G$ decomposes as a direct sum $\oplus T(p)$ and each of the
components contains either a copy of $\Z/p^k$ or one of $\Zpinf$
as a direct summand. This implies by Lemma~\ref{directsum} that
either $M(\Z/p^k, 1)$ or $M(\Zpinf, 1)$ is $M$-cellular. We
conclude now by Theorem \ref{finitetorsioncellular} and
Proposition \ref{pinfcellular}. \hfill{\qed}
\end{pf}


\section{Counter-examples to the naive characterization}
\label{sec counter}
In this section we show that the naive description of cellular
spaces fails, even in the simply connected case. It is therefore
necessary in general to work with the more complicated
characterization of Theorem~\ref{general-description} and the
notion of quasi-radical groups. Our counterexample relies on the
following computation.

\subsection{The cellularization of $K(A, 2)$}
We compute here the $M(G,1)$-cellularization of an Eilenberg-Mac
Lane space $K(A,2)$, where $A$ is any abelian group. For
simplicity we shall suppose that $K(A, 2)$ is $HR$-acyclic. This
is no longer a restriction, since for calculating $CW_M K(A,2)$ we
can first calculate $\ov P_M K(A,2)$, which is a certain product
$K(A',2)\times K(N',1)$ by  \cite[Proposition~6.1]{RS99} and then
apply $CW_M$. Note that $K(A',2)$ is $HR$-acyclic and it is easy
to calculate $CW_MK(N',1)$.

\begin{theorem}
\label{effectKA2}
Let $M=M(G,1)$ be a two-dimensional Moore space and $R$ be the
ring associated to~$G$. Suppose that $K(A,2)$ is $HR$-acyclic.
Then,
$$CW_M K(A,2)\simeq K(\Ker \varphi,2) \times K(\coker \varphi, 1)$$
where $\varphi: A \to E \to E/T_GE$, and $E$ is the ``universal"
extension for $A$ as in Definition~\ref{quasiradical}.
\end{theorem}

\begin{pf}
Let us consider the map $K(A, 2) \rightarrow K(E, 2)$. It
satisfies the conditions of Theorem~\ref{Chacholski} since its
homotopy fiber is $K(\oplus \Gab, 1)$, an $M$-cellular space by
Corollary~\ref{corollary_Gradical}. We must therefore compute
$P_{\Sigma M} K(E, 2)$ to identify the cellularization of $K(A,
2)$.

We use now Bousfield's formulas in \cite[Theorem~7.5]{Bou97}. When
$\Gab$ is a torsion group $P_{\Sigma M} K(E, 2) \simeq K(E/T_G E,
2)$. When $\Gab$ is not a torsion group, then $P_{\Sigma M} K(E,
2) \simeq K(E/T_G E, 2) \times K(B, 3)$ where $B = \prod_J
\Hom(\Z_{p^\infty}, E)$. But $K(E, 2)$ is $HR$-acyclic, $R =
\oplus_J \Z/p$ here, so $B=0$. Hence, in both cases, $CW_M K(A,
2)$ is the homotopy fiber of the map $K(A, 2) \rightarrow K(E/T_G
E, 2)$. \hfill{\qed}
\end{pf}

\subsection{Failure of the naive characterization}
Let $G=\Zup * \Z/p$ and $M=M(\Zup,1)\vee M(\Z/p,1)$. Then
$\Gab=\Zup \oplus \Z/p$ is not torsion, the associated set of
primes $J$ as introduced in Notation~\ref{notation} is the empty
set, $H=\Z$, and $R=0$, so that all spaces are $HR$-acyclic. If
the naive characterization given in the introduction was true, all
simply connected spaces would be $M$-cellular. But this is false:

\begin{theorem}
\label{theorem counterexample}
Let $G = \Zup * \Z/p$, and $M$ be the Moore space $M(\Z[1/p], 1)
\vee M(\Z/p, 1)$, so the associated ring $R$ is $0$. The space
$K(\Z, 2)$ is then $HR$-acyclic, its fundamental group is $\Z[1/p]
* \Z/p$-cellular, but $K(\Z, 2)$ is not $M$-cellular.
\end{theorem}

\begin{pf}
Let us compute explicitly $CW_M K(\Z, 2)$. We consider the
composite map $\alpha: \Z \xrightarrow{p} \Z \rightarrow
\Z^\wedge_p$, or equivalently $\Z \rightarrow \Z^\wedge_p
\xrightarrow{p} \Z^\wedge_p$, where both morphisms to the $p$-adic
integers are the completion maps. The short exact sequence of
cokernels shows that the cokernel of the composite is
$\Z^\wedge_p/\Z \times \Z/p$. Consider the fibration
$$
K(\Z^\wedge_p/\Z \times \Z/p, 1) \longrightarrow K(\Z, 2)
\xrightarrow{\alpha} K(\Z^\wedge_p, 2)\, .
$$
Since $K(\Z^\wedge_p/\Z, 1) \simeq CW_{M(\Z[1/p], 1)} K(\Z, 2)$ by
Proposition~\ref{divisiblecellularization}, it is $M$-cellular,
and so is $K(\Z/p, 1) \simeq CW_{M(\Z/p, 1)} K(\Z, 2)$. The map
$\alpha$ satisfies thus the conditions of
Theorem~\ref{Chacholski}. But notice now that the space
$K(\Z^\wedge_p, 2)$ is $\Sigma M$-null, which exhibits $CW_M K(\Z,
2)$ as $K(\Z^\wedge_p/\Z \times \Z/p, 1)$. In particular $K(\Z,
2)$ is not $M$-cellular. \hfill{\qed}
\end{pf}

The Moore space in the above counterexample has non-abelian
fundamental group. We do not know whether the naive
characterization holds in the abelian case.

\begin{question}\label{abeliancharacterization}
{\rm Let $M = M(G, 1)$ be a two-dimensional Moore space with $G$
abelian. Is it true that a space $X$ is $M$-cellular if and only
if $\pi_1 X$ is $G$-cellular and $X$ is $HR$-acyclic? The question
is even open for $G = \Zpinf$.}
\end{question}


\end{document}